\documentclass[12pt]{article}

\usepackage[T1]{fontenc}
\usepackage{amsfonts,amsthm}
\usepackage{amssymb,amsmath}

\begin{document}

\title{On Bohr-Sommerfeld-Heisenberg Quantization\thanks{
This paper is based on the lectures delivered by the second author (J.\'{S})
at the Workshop \textquotedblleft Q-days in Barcelona\textquotedblright ,
CRM Bellaterra, October 16--18, 2013.}}
\author{Richard Cushman and J\k{e}drzej \'{S}niatycki \\
Department of Mathematics and Statistics\\
University of Calgary,\\
Calgary, Alberta, Canada.}
\maketitle

\begin{abstract}
This paper presents the theory of Bohr-Sommerfeld-Heisenberg quantization of
a completely integrable Hamiltonian system in the context of geometric
quantization. The theory is illustrated with several examples.
\end{abstract}

\section{Introduction}

Most texts on quantum mechanics have a short section on the old quantum
theory. They discuss Bohr's quantization of the harmonic oscillator and
Sommerfeld's results on the energy spectrum of the hydrogen atom. Usually
they mention of Heisenberg's quantum mechanics and give a description of 
Schr\"{o}dinger's wave mechanics. Schr\"{o}dinger's theory is further discussed
in the framework of the modern quantum mechanics. Heisenberg's theory is
relegated to a criptic remark that Dirac proved that that the theories of
Heisenberg and of Schr\"{o}dinger are equivalent. In \cite{dirac2}, Dirac
showed that Heisenberg's matrices can be also obtained in the 
Schr\"{o}dinger theory, but he did not state that these theories give the same
physical results.

Geometric quantization provides an explanation of Dirac's theory in the
framework of modern differential geometry. Within geometric quantization, it
is easy to understand Bohr-Sommerfeld quantization rules; I discussed them
in my book \cite{sniatycki80}. However, for a long time I could not fit
Heisenberg's matrix mechanics into the framework of geometric quantization.

The breakthrough came when Richard Cushman explained to me the notion of
quantum monodromy introduced in a joint paper with Duistermaat \cite{cushman-duistermaat}. 
If a Hamiltonian system with $n$-degrees of freedom
admits a global Hamiltonian action of the $n$-torus, then the
Bohr-Sommerfeld conditions are global and they define the structure of an 
$n$-dimensional lattice on the corresponding basis of the space of quantum
states. Cushman and Duistermaat showed that, in the presence of classical
monodromy, this lattice structure was only local.

In this lecture, I will describe our understanding of Heisenberg's quantum
mechanics within the framework of geometric quantization. I do not know if
our approach has any relation to Heisenberg's ideas. However, I hope to
convince you that we obtain a well defined quantum theory consistent with
the principles of geometric quantization. More precisely, the theory we
obtain generalizes geometric quantization, as formulated by Kostant, to the
case of a singular polarization.

\section{Completely integrable systems}

Let $(P,\omega )$ be a symplectic manifold of dimension $2n.$ We consider a
completely integrable system on $(P,\omega )$ with action angle coordinates 
$(A_{i},\varphi _{i})$ defined on an open dense subset $U$ of $P$. The
symplectic form $\omega $ restricted to $U$ is $\omega _{\mid U}=d\theta $,
where $\theta =\sum_{i=1}^{n}d(A_{i}d\varphi _{i})$.

\begin{description}
\item[Assumption 1] We assume that the action coordinates $A_{i}$ are
globally \linebreak defined on $P$.
\end{description}

This implies that we have a symplectic action of the the torus group $\mathbb{T}^{n}$ 
with the momentum map $J:P\rightarrow 
\mathbb{R} ^{n}:p\mapsto J(p)= (A_{1}(p),...,A_{n}(p))$, where we have identified the
Lie algebra of $\mathbb{T}^{n}$ with $\mathbb{R}^{n}$.

\section{Bohr-Sommerfeld Quantization}

The Hamiltonian vector field $X_{f}$ of a function $f\in C^{\infty }(P)$ is
defined by $X_{f}{\mbox{$ \rule {5pt} {.5pt}\rule {.5pt} {6pt} \, $}} \omega
=-df$, where {\mbox{$ \rule {5pt} {.5pt}\rule {.5pt} {6pt} \, $}} is the
left interior product (contraction on the left).

For each $i=1,...,n$, the Hamiltonian vector field $X_{A_{i}}$ generates the
action on $P$ of the $i^{\mathrm{th}}$ component $\mathbb{T}_{i}$ of the
torus group $\mathbb{T}^{n}=\mathbb{T}\times \mathbb{T}\times ...\times 
\mathbb{T}$. We denote by $O_{i,p}$ the orbit of $\mathbb{T}_{i}$ through $p\in P.$ Clearly, $A_{i}$ is constant on each orbit $O_{i,p}.$

\begin{description}
\item[Bohr-Sommerfeld Quantization Rule] For each $i=1,...,n,$ the quantum
spectrum of $A_{i}$ consists of the values $A_{i}(p)$ on orbits $O_{i,p}$
satisfying the condition 
\begin{equation}
\int_{O_{i,p}}A_{i}d\varphi _{i}=m_{i}h,  \label{eq-one}
\end{equation}
where $m_{i}$ is the integer and $h$ denotes Planck's constant.
\end{description}

Integrating equation (\ref{eq-one}), we conclude that the quantum spectrum
of the $i^{\mathrm{th}}$ action is given by 
\begin{equation}
A_{i}=m_{i}\hbar ,  \label{eq-twoa}
\end{equation}
where $\hbar $ is Planck's constant divided by $2\pi .$

\begin{description}
\item[Assumption 2] For each n-tuple $\mathbf{m}=(m_{1},...,m_{n})$ of
integers, the set 
\begin{equation}
\mathbb{T}_{\mathbf{m}}=\{p\in P\mid A_{i}(p)=m_{i}\hbar \text{~~}\forall
~~i=1,...,n\}.  \label{eq-three}
\end{equation}
is connected.
\end{description}

Under this assumption, $\mathbb{T}_{\mathbf{m}}$ is a torus. Otherwise, it
would be the union of disjoint tori, and we would have to introduce an
additional index to label connected components. In the following, we shall
refer to sets $\mathbb{T}_{\mathbf{m}}$ defined by equation (\ref{eq-three})
as Bohr-Sommerfeld tori.

\section{Link to Geometric Quantization}

Suppose that we want to perform geometric quantization of our completely
integrable system in the real polarization $D$ spanned by the Hamiltonian
vector fields $X_{A_{i}}$ of the momenta $A_{1},...,A_{n}$.

Let $L$ be a prequantization line bundle of $(P,\omega )$. Thus, $L$ is a
complex line bundle over $P$, with a connection $\nabla $ such that 
\begin{equation*}
(\nabla _{X_{1}}\nabla _{X_{2}}-\nabla _{X_{2}}\nabla _{X_{1}}-\nabla
_{\lbrack X_{1},X_{2}]})\sigma =-\frac{i}{\hbar }\omega (X_{1},X_{2}) 
\end{equation*}%
for each section $\sigma $ of $L$ and every pair $X_{1},X_{2}$ of vector
fields on $P.$

The quantum states of the system are given by sections $\sigma $ of $L$ that are
covariantly constant along the polarization $D.$ If $\Lambda $ is a leaf of $D$, it is a torus, and the restriction $\sigma _{\mid \Lambda }$ of a
section $\sigma $ of $L$ that is covariantly constant along $D$ vanishes
unless the holonomy group of the restriction of $\nabla $ to $\Lambda $
vanishes.

\begin{description}
\item[Proposition] The holonomy group of the restriction of $\nabla $ to a
leaf $\Lambda $ of $D$ vanishes if and only if $\Lambda $ satisfies the
Bohr-Sommerfeld conditions; that is $\Lambda =\mathbb{T}_{\mathbf{m}}$ for
some $n$-tuple $\mathbf{m}=(m_{1},...,m_{n})\in \mathbb{Z}^{n}$.
\end{description}

The proof of this proposition can be found in reference \cite{sniatycki80}.
\hfill $\square $ \medskip

For a completely integrable system, tori $\mathbb{T}_{\mathbf{m}}$ are
submanifolds of $P$ of codimension at least $n$, and the only smooth section
of $L$ that is covariantly constant along $D$ is the zero section. However,
we may identify quantum states with distribution sections of $L$ that are
smooth and covariantly constant along leaves of $D$. Under this
interpretation of quantum states, to every non-empty Bohr-Sommerfdeld torus 
$\mathbb{T}_{\mathbf{m}}$ in $P$, we can associate a non-zero distribution
section $\sigma _{\mathbf{m}}$ of $L$ with support in 
$\mathbb{T}_{\mathbf{m}},$ and such that the restriction 
$\sigma _{\mathbf{m}\mid \mathbb{T}_{\mathbf{m}}}$ of $\sigma _{\mathbf{m}} $ to 
$\mathbb{T}_{\mathbf{m}}$ is a
smooth covariantly constant section the restriction of $L$ to $\mathbb{T}_{\mathbf{m}}$. On the space $\mathfrak{S}$ of distribution sections of $L$
that are spanned by sections $\sigma _{\mathbf{m}}$ introduce a scalar
product $\left\langle \cdot \mid \cdot \right\rangle $ such that 
\begin{equation}
\left\langle \sigma _{\mathbf{m}}\mid \sigma _{\mathbf{m}^{\prime}}\right\rangle = \delta _{\mathbf{mm}^{\prime }}=\delta _{m_{1}m_{1}^{\prime
}}....\delta _{m_{n}m_{n}^{\prime }}.  \label{eq-four}
\end{equation}

For each Bohr-Sommerfeld torus $\mathbb{T}_{\mathbf{m}}$, the section 
$\sigma _{\mathbf{m}}$ introduced above is defined by $\mathbb{T}_{\mathbf{m}}
$ up to an arbitrary non-zero complex factor. Therefore, the collection 
$\{\mathbb{T}_{\mathbf{m}}\}$ of all Bohr-Sommerfeld tori in $P$ determines
only the orthogonality property of basic vectors $\sigma _{\mathbf{m}}$. For
a covariantly constant section section $\sigma $ with suppoort in 
${\mathbb{T}}_{\mathbf{m}}$, the norm $\left\Vert \sigma \right\Vert $ depends on the
choice of the basic section $\sigma _{\mathbf{m}}$.

We denote by $\mathfrak{H}$ the Hilbert space obtained by the completion of $\mathfrak{S}$ in the norm given by $\left\langle \cdot \mid \cdot
\right\rangle $. $\mathfrak{H}$ is our space of quantum states. To each
function $f\in C^{\infty }(P)$, such that $f=F(A_{1},...,A_{n})$ for some 
$F\in C^{\infty }(\mathbb{R}^{n})$, the Bohr-Sommerfeld quantization
associates the quantum operator $\mathbf{Q}_{f}$ on $\mathfrak{H} $ such
that, for every basic section $\sigma _{\mathbf{m}}$, 
\begin{equation*}
\mathbf{Q}_{f}\sigma _{\mathbf{m}}=F(\mathbf{m}\hbar )\sigma _{\mathbf{m}}. 
\end{equation*}%
It follows from Assumption 2 that the spectrum of the action operators $\mathbf{Q}_{A_{i}}$ is simple.

A shortcoming of the Bohr-Sommerfeld quantization is that it is defined only
on the commutative algebra consisting of smooth smooth functions of the
actions. In particular, Bohr-Sommerfeld quantization does not allow for
quantization of any function of the angles. Moreover, it leads only to
diagonal operators in $\mathfrak{H}$.

\section{Shifting operators}

Bohr-Sommerfeld conditions together with Assumptions 1 and 2 imply that the
basis $\{\sigma _{\mathbf{m}}\}$ is a lattice. Therefore, there are well
defined operators corresponding to shifting along the generators of the
lattice.

For each $i=1,...,n$, let 
\begin{equation*}
\mathbf{m}_{i}=\{m_{1},...,m_{i-1},m_{i}-1,m_{i+1},...,m_{n}\} 
\end{equation*}
and 
\begin{equation*}
\mathbf{m}^{i}=\{m_{1},...,m_{i-1},m_{i}+1,m_{i+1},...,m_{n}\}. 
\end{equation*}
We define shifting operators $\mathbf{a}_{i}$ on $\mathfrak{H}$ by 
\begin{equation}
\mathbf{a}_{i}\sigma _{\mathbf{m}}=\left\{ 
\begin{array}{ccc}
\sigma _{\mathbf{m}_{i}} & if & \mathbb{T}_{\mathbf{m}_{i}}\neq \emptyset \\ 
0 & if & \mathbb{T}_{\mathbf{m}_{i}}=\emptyset
\end{array}
\right. .  \label{eq-five}
\end{equation}
The adjoint operators $\mathbf{a}_{i}^{\dagger }$ are given by 
\begin{equation}
\mathbf{a}_{i}^{\dagger }\sigma _{\mathbf{m}}=\left\{ 
\begin{array}{ccc}
\sigma _{\mathbf{m}^{i}} & if & \mathbb{T}_{\mathbf{m}^{i}}\neq \emptyset \\ 
0 & if & \mathbb{T}_{\mathbf{m}^{i}}=\emptyset
\end{array}
\right. .  \label{eq-six}
\end{equation}

\begin{description}
\item[Proposition 1] The shifting operators satisfy the following
commutation \linebreak relations
\begin{subequations}
\begin{eqnarray}
\lbrack \mathbf{a}_{k},\mathbf{Q}_{A_{j}}] &=&\delta _{kj}\hbar 
{\mathbf{a}}_{k},  \label{eq-sevena} \\
\lbrack \mathbf{a}_{k}^{\dagger },\mathbf{Q}_{A_{j}}] &=&-\delta _{kj}\hbar 
\mathbf{a}_{k}^{\dagger }.  \label{eq-sevenb}
\end{eqnarray}
\end{subequations}
\end{description}
The Poisson bracket relations between actions and angles are 
\begin{equation*}
\{e^{-i\varphi _{k}},A_{j}\}=-i\delta _{kj}e^{-i\varphi _{k}}. 
\end{equation*}
Hence, Dirac's quantization conditions 
\begin{equation}
\lbrack \mathbf{Q}_{f_{1}},\mathbf{Q}_{f_{2}}]=i\hbar \mathbf{Q}
_{\{f_{1},f_{2}\}}  \label{9}
\end{equation}
suggest the identification $\mathbf{a}_{k}=\mathbf{Q}_{e^{-i\varphi _{k}}}$
and $\mathbf{a}_{k}^{\dagger }= \mathbf{Q}_{e^{i\varphi _{k}}}$, where $\varphi _{k}$ is the angle coordinate corresponding to the action $A_{k}$,
provided the functions $e^{-i\varphi _{k}}$ and $e^{i\varphi _{k}}$ are
globally defined on $P$.

\section{Heisenberg Quantization}

Since not all sets $\mathbb{T}_{\mathbf{m}}$ are $n$-tori, we cannot expect
that all exponential functions $e^{-i\varphi _{k}}$ are globally defined. We
can try to replace $e^{-i\varphi _{k}}$ by a globally defined smooth
function $f_{k}$ of the form $\chi _{k}=r_{k}e^{-i\varphi _{k}}$, where the
coefficient $r_{k}$ depends only on the actions and vanishes at the points
at which $e^{i\varphi _{k}}$ is not defined. In the following we shall refer
to functions $\chi _{k}$ as Heisenberg functions.

We have the following Poisson bracket relations
\begin{equation}
\{\chi _{k},A_{j}\}=-i\delta _{kj}\chi _{k}\text{ and }\{{\bar{\chi}}_{k},A_{j}\}=i\delta _{kj}\bar{\chi}_{k}.  \label{10a}
\end{equation}
By Dirac's quantization conditions, we get 
\begin{eqnarray}
\lbrack \mathbf{Q}_{\chi _{k}},\mathbf{Q}_{A_{j}}] &=&\delta _{kj}\hbar 
\mathbf{Q}_{\chi _{k}},  \label{11} \\
\lbrack \mathbf{Q}_{\bar{\chi}_{k}},\mathbf{Q}_{A_{j}}] &=&-\delta
_{kj}\hbar \mathbf{Q}_{\bar{\chi}_{k}}.  \label{12}
\end{eqnarray}
For each basic vector $\sigma _{\mathbf{m}}$ of $\mathfrak{H},$ 
\begin{eqnarray}
\mathbf{Q}_{A_{j}}(\mathbf{Q}_{\chi _{j}}\sigma _{\mathbf{m}}) &=&
{\mathbf{Q}}_{\chi _{j}}(\mathbf{Q}_{A_{j}}\sigma _{\mathbf{m}})-[\mathbf{Q}_{\chi _{j}},
\mathbf{Q}_{A_{j}}]\sigma _{\mathbf{m}}  \label{12a} \\
&=&\mathbf{Q}_{\chi _{j}}(\hbar m_{j}\sigma _{\mathbf{m}})-\hbar {\mathbf{Q}}_{\chi _{j}}
\sigma _{\mathbf{m}}  \notag \\
&=&\hbar (m_{j}-1)\mathbf{Q}_{\chi _{j}}\sigma _{\mathbf{m}}.  \notag
\end{eqnarray}
Thus, $\mathbf{Q}_{\chi _{j}}\sigma _{\mathbf{m}}$ is proportional to 
$\sigma _{\mathbf{m}_{j}}.$ A similar argument shows that 
$\mathbf{Q}_{\bar{\chi}_{j}}\sigma _{\mathbf{m}}$ is proportional to \thinspace 
$\sigma _{\mathbf{m}^{j}}.$ Hence, $\mathbf{Q}_{\chi _{j}}$ and 
$\mathbf{Q}_{\bar{\chi}_{j}}$ act as shifting operators, namely, 
\begin{equation}
\mathbf{Q}_{\chi _{j}}\sigma _{\mathbf{m}}=b_{\mathbf{m},j}
\sigma _{\mathbf{m}_{j}}\text{ and }\mathbf{Q}_{\bar{\chi}_{j}}\sigma _{\mathbf{m}}=
c_{\mathbf{m},j}\sigma _{\mathbf{m}^{j}}  \label{12b}
\end{equation}
for some coefficients $b_{\mathbf{m},j}$ and $c_{\mathbf{m},j}$.

We can use Dirac's quantization conditions 
\begin{equation}
\lbrack \mathbf{Q}_{\chi _{j}},\mathbf{Q}_{\chi _{k}}]=i\hbar 
\mathbf{Q} _{\{\chi _{j},\chi _{k}\}}\text{ ~~and~~ }[\mathbf{Q}_{\chi _{j}},\mathbf{Q}%
_{\bar{\chi}_{k}}]=i\hbar \mathbf{Q}_{\{\chi _{j},\bar{\chi}_{k}\}}
\label{12c}
\end{equation}
and the identification 
\begin{equation}
\mathbf{Q}_{\chi _{j}}^{\dagger }=\mathbf{Q}_{\bar{\chi}_{j}}  \label{12d}
\end{equation}
to determine the coefficients $b_{\mathbf{m},j}$ and $c_{\mathbf{m},j}$,
which must satisfy the consistency conditions:
\begin{equation}
b_{\mathbf{m},j}=0\text{ if }\mathbb{T}_{\mathbf{m}_{j}}=\emptyset \text{ ~\
and ~\ }c_{\mathbf{m},j}=0\text{ if }\mathbb{T}_{\mathbf{m}^{j}}=\emptyset .
\label{12e}
\end{equation}

The Bohr-Sommerfeld-Heisenberg quantization described here is an extension
of the Bohr-Sommerfeld theory. In the Bohr-Sommerfeld-Hesenberg
quantization, the Hilbert space $\mathfrak{H}$ of quantum states is the same
as in the Bohr-Sommerfeld theory. However, in the Bohr-Sommerfeld-Heisenberg
theory, we can quantize functions that first degree polynomials in $\chi
_{k} $ and $\bar{\chi}_{k}$ with coefficients given by smooth functions of
the actions:
\begin{equation*}
F(A_{1},...,A_{n})+\sum_{k=1}^{n}[F_{k}(A_{1},....,A_{n})\chi _{k}+\tilde{F} _{k}(A_{1},....,A_{n})\bar{\chi}_{k}]. 
\end{equation*}
The resulting operators on $\mathfrak{H}$ first degree polynomials in
shifting operators. Higher powers of shifting operators are well defined on 
$\mathfrak{H}$, but they need not be quantizations of the corresponding
powers of the functions $f_{k}$ or $\bar{f}_{k}$ (the usual factor ordering
problem).

\section{Examples}

\subsection{1-dimensional harmonic oscillator}

The phase space of the 1-dimensional harmonic oscillator is $P=\mathbb{R}^{2}
$ with coordinates $(p,q)$ and the symplectic form $\omega =dp\wedge dq$.
The Hamiltonian is $H=\frac{\scriptstyle1}{\scriptstyle2}(p^{2}+q^{2})$. In
polar coordinates $(p,q)=(r\cos \varphi ,r\sin \varphi )$, where $r=
\sqrt{p^{2}+q^{2}}$ and $\varphi =\tan \frac{q}{p}$, we have $\omega =dH\wedge
d\varphi $. Here $H=\frac{1}{2}r^{2}$ is the action variable, while $\varphi 
$ is the corresponding angle. The Heisenberg function $\chi
=p-iq=re^{-i\varphi }$ leads to quantization equivalent to the Bargmann
quantization \cite{bargmann}. It should be noted that $r=\sqrt{2H}$ is not a
smooth function of $H$, but $\chi $ is in $C^{\infty }(P)$. For full
details see \cite{cushman-sniatycki13}. 

\subsection{ Coadjoint orbits of SO(3)}

Following Souriau \cite{souriau} we use the presentation of coadjoint orbits
of $\mathrm{SO}\,(3)$ spheres $S_{r}^{2}=\{(x^{1},x^{2},x^{3})\in 
{\mathbb{R}}^{3}\mid (x^{1})^{2}+(x^{2})^{2}+(x^{3})^{2}=r^{2}\}$ endowed with a
symplectic form $\omega =\frac{\scriptstyle1}{\scriptstyle r}{\mathrm{vol}}_{S_{r}^{2}}$, where 
$\mathrm{vol}_{S_{r}^{2}}$ is the standard area form on 
$S_{r}^{2}$ with $\int_{S_{r}^{2}}\mathrm{vol}_{S_{r}^{2}}=4\pi r^{2}$. A
coadjoint orbit $S_{r}^{2}$ is qantizable if $r=\frac{n}{2}\hbar $, where $n$
is an integer.

For each $i=1,2,3,$ we denote by $J^{i}$ the restriction of $x^{i}$ to the
sphere $S_{r}^{2}$. The functions $J^{1},$ $J^{2}$ and $J^{3}$ are
components of the momentum map of the co-adjoint action. They satisfy the
Poisson bracket relations $\{J^{i},J^{j}\}=\sum_{k=1}^{3}\varepsilon
_{ijk}J^{k}$. In spherical polar coordinates 
\begin{equation*}
J^{1}=r\sin \theta \cos \varphi \text{, \ }J^{2}=r\sin \theta \sin \varphi 
\text{, \ }J^{3}=r\cos \theta ,
\end{equation*}
and 
\begin{equation*}
\omega =r\sin \theta d\varphi \wedge d\theta =-d(r\cos \theta d\varphi
)=d(J^{3}d(-\varphi )).
\end{equation*}
Thus, $(J^{3},-\varphi )$ are action-angle coordinates for an integrable
system $(J^{3},S_{r}^{2},\omega )$. In this case, a Heisenberg function is 
$\chi =J_{+}=\sqrt{r^{2}-(J^{3})^{2}}e^{i\varphi }$, and the resulting
Bohr-Sommerfeld-Heisenberg quantization leads to the irreducible unitary
representation of $SO(3)$ corresponding to the co-adjoint orbit $S_{r}^{2}$.
For more details, see \cite{cushman-sniatycki13}. The presented treatment
closely resembles the approach of Schwinger \cite{schwinger}.

\subsection{2-dimensional harmonic oscillator}

The configuration space of the $2$-dimensional harmonic oscillator is $\mathbb{R}^{2}$ with coordinates $x=(x^{1},x^{2})$. The phase space is $T^{\ast }\mathbb{R}^{2}=\mathbb{R}^{4}$ with coordinates $(x,y)=(x^{1},x^{2},y^{1},y^{2})$ and canonical symplectic form $\omega
=d(y^{1}dx^{1}+y^{2}dx^{2})$. The Hamiltonian function of the $2$-dimensional harmonic oscillator is 
\begin{equation*}
H(x,y)=\frac{\scriptstyle 1}{\scriptstyle 2}\big( (x^{1})^{2}+(x^{2})^{2}
\big) + \frac{\scriptstyle 1}{\scriptstyle 2}\big( (y^{1})^{2}+(y^{2})^{2}
\big) . 
\end{equation*}
Orbits of the Hamiltonian vector field $X_{H}$ of $H$ are periodic of period 
$2\pi $. The function $L(x,y)=x^{1}y^{2}-x^{2}y^{1}$ generates an action of $S^{1}$ on 
$T^{\ast }\mathbb{R}^{2}$ that preserves the Hamiltonian $H$.
Hence, $(H,L,T^{\ast }\mathbb{R}^{2},\omega )$ is a completely integrable
system. Let 
\begin{equation}
\begin{array}{lcl}
x_1 = \mbox{$\frac{{\scriptstyle -1}}{{\scriptstyle \sqrt{2}}}$}(r_1 
\cos {\vartheta }_1 + r_2\cos {\vartheta }_2) & \, \, \, & y_1 = 
\mbox{$\frac{{\scriptstyle 1}}{{\scriptstyle \sqrt{2}}}$}(r_1 
\sin {\vartheta }_1 + r_2\sin {\vartheta }_2) \\ 
\rule{0pt}{16pt} x_2 = 
\mbox{$\frac{{\scriptstyle 1}}{{\scriptstyle
\sqrt{2}}}$}(-r_1 \sin {\vartheta }_1 + r_2\sin {\vartheta }_2) & \, \, \, & 
y_2 = \mbox{$\frac{{\scriptstyle 1}}{{\scriptstyle \sqrt{2}}}$}(-r_1 \cos {\vartheta }_1 + 
r_2\cos {\vartheta }_2) .
\end{array}
\label{eq-emactang}
\end{equation}
be a change of coordinates from rectangular $(x,y)$ variables to polar
variables $(r_1, r_2, {\vartheta }_1, {\vartheta }_2)$. A computation shows
that $H(r,\vartheta ) = \mbox{$\frac{\scriptstyle
1}{\scriptstyle 2}\,$} (r^2_1 +r^2_2)$ and $L(r,\vartheta ) = 
\mbox{$\frac{\scriptstyle 1}{\scriptstyle 2}\,$} (r^2_1-r^2_2)$ and that the
change of coordinates (\ref{eq-emactang}) pulls back the symplectic form 
$\omega = d y_1 \wedge d x_1 + dy_2 \wedge d x_2$ to the symplectic form $\Omega = 
d (\mbox{$\frac{\scriptstyle 1}{\scriptstyle 2}\,$} r^2_1) \wedge d 
{\vartheta }_1 + d (\mbox{$\frac{\scriptstyle 1}{\scriptstyle 2}\,$} r^2_2)
\wedge d{\vartheta }_2$. Let $A_1 = \mbox{$\frac{\scriptstyle
1}{\scriptstyle 2}\,$} r^2_1 = \mbox{$\frac{\scriptstyle 1}{\scriptstyle
2}\,$} \big( E(r, \vartheta ) + L(r, \vartheta ) \big) \ge 0 $ and $A_2 = 
\mbox{$\frac{\scriptstyle 1}{\scriptstyle 2}\,$} r^2_2 = 
\mbox{$\frac{\scriptstyle 1}{\scriptstyle 2}\,$} \big( E(r, \vartheta ) - L(r, \vartheta
) \big) \ge 0$. Then $(A_1, A_2, {\vartheta }_1, {\vartheta }_2)$ with $A_1
>0$, and $A_2 >0$ and symplectic form $\Omega = d A_1\wedge d {\vartheta }_1
+ d A_2 \wedge d {\vartheta }_2$ are real analytic action-angle coordinates
for the $2$-dimensional harmonic oscillator. These \linebreak 
coordinates extend real
analytically to the closed domain $A_1 \ge 0$ and $A_2 \ge 0$. The
Heisenberg functions $\chi _{1}=r_1{\mathrm{e}}^{i{\vartheta}_1}$ and $\chi
_{2}=r_2{\mathrm{e}}^{i{\vartheta}_2}$ give rise to the
Bohr-Sommerfeld-Heisenberg quantization of the $2$-dimensional harmonic
oscillator. For more details see \cite{cushman-sniatycki12}.

\subsection{Mathematical Pendulum}

The phase space of the mathematical pendulum is $T^{\ast }S^{1}$ with
coordinates $(p,\alpha )$ and symplectic form $\omega =dp\wedge d\alpha $.
The Hamiltonian of the system is $H=\frac{\scriptstyle1}{\scriptstyle2}
p^{2}-\cos \alpha +1$. The Hamiltonian system $(H,T^{\ast }S^{1},\omega )$
violates \linebreak 
Assumption 2, because for $H>2,$ level sets of the
Hamiltonian $H$ have two connected components. We are investigating how to
extend to this case the theory presented here.

\end{document}